# Numerical Solution of Fuzzy Stochastic Differential Equation

S. Nayak and S. Chakraverty[1]

Department of Mathematics, National Institute of Technology, Rourkela, Odisha -769008, India

**Abstract**

In this paper an alternative approach to solve uncertain Stochastic Differential Equation (SDE) is proposed. This uncertainty occurs due to the involved parameters in system and these are considered as Triangular Fuzzy Numbers (TFN). Here the proposed fuzzy arithmetic in [2] is used as a tool to handle Fuzzy Stochastic Differential Equation (FSDE). In particular, a system of Ito stochastic differential equations is analysed with fuzzy parameters. Further exact and Euler Maruyama approximation methods with fuzzy values are demonstrated and solved some standard SDE.

1. **Introduction**

The concept of Stochastic Differential Equation (SDE) has been initiated by Einstein in 1905 [15]. In his article he presented a mathematical connection between microscopic random motion of particles and the macroscopic diffusion equation. Later it has been seen that the stochastic differential equation (SDE) model plays a prominent role in a range of application areas such as physics, chemistry, mechanics, biology, microelectronics, economics and finance. Earlier the SDEs were solved by using Ito integral as an exact method which is discussed in [9]. But using exact method it is noticed that there occur some difficulty to study nontrivial problems and hence approximation methods are used. In this context various authors have given their contribution in these field but we have mentioned which are directly related to this problem. In 1982, Rumelin [14] defined general Runge-Kutta approximations for the solution of stochastic differential equations and there was given an explicit form of the correction term. This work was carried out and then Kloeden and Platen [8] are discussed about the numerical solutions of stochastic differential equation in detail. Platen [13] added discrete time strong and weak approximation methods for the numerical methods to find the solution of stochastic differential equations. Next, Higham [5] gave a major contribution in this field to solve the approximate solutions of stochastic differential equations and discussed few problems. Further Higham and Kloeden [6] investigated nonlinear stochastic differential equations numerically. They presented two implicit methods for Ito stochastic differential equations (SDEs) with Poisson-driven jumps. The first method is a split-step extension of the backward Euler method and the second method arises from the introduction of a compensated, martingale, form of the Poisson process. In this context different authors have tried for various other diffusion and application based problems. Hayes and Allen [4] solved stochastic point kinetic reactor problem. They modelled the point stochastic reactor problem into ordinary time dependent stochastic differential equation and studied the stochastic behaviour of the neutron flux.

---
[1] corresponding author
E-mail address: sne_chak@yahoo.com (S. Chakraverty), sukantgacr@gmail.com (S. Nayak)

It may be noted from the literature review that, authors have discussed the stochastic differential equations which contain crisp parameters. But in general the involved parameters may not be crisp rather these may be uncertain. Here the uncertain parameters are considered as TFN. In this context, Kim [7] considered fuzzy sets space for real line and the existence and uniqueness of the solution is obtained. The solution is investigated by taking particular conditions which are imposed on the structure of integrated fuzzy stochastic processes such that a maximal inequality for fuzzy stochastic Ito integral holds. Next, Ogura [11] proposed an approach to solve FSDE which does not contain any notion of fuzzy stochastic Ito integral and the method was based on the selections of sets. Further, Malinowski and Michta [9] presented the existence and uniqueness of solutions to the FSDEs driven by Brownian motion and the continuous dependence on initial condition and stability properties are established.

In this investigation, a general approach has been described to handle fuzzy numbers associated with the FSDE. The concept of fuzzy stochastic Ito integral has also been used to obtain the solutions of exact method. Also, using fuzzy arithmetic [2, 3, 10] for TFNs, numerical solutions of standard FSDEs are investigated.

For the sake of completeness initially we have discussed the crisp SDE and it is solved analytically through Ito integral techniques. It is also been noted that there was difficulty to handle nontrivial problems using analytical method so we have used numerical method to solve. Further the same problems are discussed for uncertain cases and corresponding FSDE are solved. The obtained results are shown graphically and the uncertain width of the solution is discussed.

## 2. Preliminary

Let us consider a standard stochastic differential equation
$$dX = a(t, X)dt + b(t, X)dW_t \qquad (1)$$
where Eq. (1) is written in differential form.
The integral form of Eq. (1) becomes
$$X(t) = X(0) + \int_0^t a(s, y)ds + \int_0^t b(s, y)dW_s \qquad (2)$$
where the last term in the right hand side of Eq. (2) is called Ito integral.

We take $c = t_0 < t_1 < t_2 < ... < t_{n-1} < t_n = d$ be a grid of points on an interval $[c, d]$, then Ito integral may be defined in the following limit form
$$\int_c^d f(t)dW_t = \lim_{\Delta t \to 0} \sum_{i=1}^{n} f(t_{i-1})\Delta W_i \qquad (3)$$
where $\Delta W_i = W_{t_i} - W_{t_{i-1}}$, a step of Brownian motion across the interval.

## 3. Analytical solution of Stochastic Differential Equations (SDE)

Let us consider the Eq. (1), which is first solved analytically by using Ito formula.
Ito formula says that if $X_t$, an Ito process given by
$$dX_t = udt + vdW_t \qquad (4)$$

Let $g(t, x) \in C^2([0, \infty] \times \Re)$ (i.e. $g$ is twice continuous differentiable on $[0, \infty] \times \Re$). Then $Y_t = g(t, X_t)$ is again Ito process and

$$dY_t = \frac{\partial g}{\partial t}(t, X_t)dt + \frac{\partial g}{\partial x}(t, X_t)dX_t + \frac{1}{2}\frac{\partial^2 g}{\partial x^2}(t, X_t)(dX_t)^2 \qquad (5)$$

where $(dX_t)^2 = (dX_t).(dX_t)$ is computed as follows

$$dt.dt = dt.dW_t = dW_t.dt = 0;$$
$$dW_t.dW_t = dt.$$

Example 1

Let us consider a stochastic differential equation [9]

$$\begin{cases} \dfrac{dX_t}{dt} = a_t X_t, \\ X(0, x) = X_0 \end{cases} \qquad (6)$$

where $a_t = r_t + \alpha W_t$

$W_t$ and $\alpha$ are noise and constant respectively.

Eq. (6) may be written as

$$\frac{dX_t}{dt} = (r_t + \alpha W_t)X_t$$
$$= r_t X_t + \alpha W_t X_t$$
$$\Rightarrow dX_t = r_t X_t dt + \alpha X_t dB_t \quad (\because W_t.dt = dB_t)$$
$$\Rightarrow \frac{dX_t}{X_t} = r_t dt + \alpha dB_t$$
$$\Rightarrow \int_0^t \frac{dX_s}{X_s} = r_t t + \alpha B_t$$

Using Ito formula for the function $g(t, x) = \ln x$ [12], we get the following

$$d(\ln X_t) = \frac{1}{X_t}dX_t - \frac{1}{2}\left(\frac{1}{X_t^2}\right)(dX_t)^2$$
$$= \frac{1}{X_t}dX_t - \frac{1}{2}\left(\frac{1}{X_t^2}\right)\alpha^2 X_t^2 dt$$
$$= \frac{1}{X_t}dX_t - \frac{1}{2}\alpha^2 dt$$

Integrating the above we get,

$$\int_0^t d(\ln X_t) = \int \frac{dX_t}{X_t} - \int \frac{1}{2}\alpha^2 dt$$

$$\Rightarrow \ln X_t - \ln X_0 = \int r_t dt + \alpha dB_t - \frac{1}{2}\alpha^2 t$$

$$\Rightarrow \ln \frac{X_t}{X_0} = \left(r_t - \frac{1}{2}\alpha^2\right)t + \alpha B_t$$

$$\Rightarrow X_t = X_0 e^{\left(r_t - \frac{1}{2}\alpha^2\right)t + \alpha B_t}$$

Here we found that except some standard problems exact method may not be applicable for others. Hence we need numerical treatment to handle non trivial problems and this is discussed in the following sections.

## 4. Numerical solution of Stochastic Differential Equations (SDE)

To find the numerical solution of the stochastic differential equation let us assign a grid of points, $c = t_0 < t_1 < t_2 < ... < t_{n-1} < t_n = d$ and approximate $x$ values $w_0 < w_1 < w_2 < ... < w_n$ to be determined at the respective $t$ points.

Consider SDE initial value problem [1]

$$\begin{cases} dX(t) = a(t,X)dt + b(t,X)dW_t \\ X(0) = X_0 \end{cases} \quad (7)$$

Eq. (7) is solved numerically as follows.

**Euler-Maruyama Method**

Let us consider a time discrete approximation of an Ito process for the stochastic differential equation (SDE)

$$\begin{cases} dX(t) = a(t,X)dt + b(t,X)dW_t \\ X(c) = X_c \end{cases} \quad (8)$$

Then the approximation Euler-Maruyama scheme may be represented in the following manner

$$\begin{aligned} X_0 &= w_0 \\ w_{i+1} &= w_i + a(t_i, w_i)\Delta t_{i+1} + b(t_i, w_i)\Delta W_i \end{aligned} \quad (9)$$

where

$$\Delta t_{i+1} = t_{i+1} - t_i$$
$$\Delta W_{i+1} = W(t_{i+1}) - W(t_i)$$

Define $N(0,1)$ be the normal distribution and each random number $\Delta W_i$ is computed as

$$\Delta W_i = z_i \sqrt{\Delta t_i}$$

where, $z_i$ is chosen from $N(0,1)$.

Here the obtained set $\{w_0, w_1, ..., w_n\}$ is an approximation realization of the solution stochastic process $X(t)$ which depends on the random numbers $z_i$ that were chosen. Since, $W_t$ is a stochastic process, each realization will be different and so will our approximations.

## 5. Fuzzy arithmetic

Fuzzy set $\tilde{A}$ is the collection of pair of elements and its membership functions. The membership function $\mu_{\tilde{A}}$ is defined as

$$\mu_{\tilde{A}} : X \to [0,1]$$

where $X$ is the universal set. Fuzzy number is a convex normal set which are described as follow.

A fuzzy number $z$ in parametric form is a pair $[\underline{z}, \bar{z}]$ of functions $\underline{z}(\alpha)$ and $\bar{z}(\alpha)$, where $\alpha \in [0,1]$ which satisfy the following conditions.

i. $\underline{z}(\alpha)$ is a bounded non-decreasing left continuous function in (0, 1], and right continuous at 0.

ii. $\bar{z}(\alpha)$ is a bounded non-increasing left continuous function in (0, 1], and right continuous at 0.
iii. $\underline{z}(\alpha) \leq \bar{z}(\alpha), 0 \leq \alpha \leq 1$.

A fuzzy number $\tilde{A} = [a^L, a^N, a^R]$ is said to be triangular fuzzy number (Fig. 1) when the membership function is given by

$$\mu_{\tilde{A}}(x) = \begin{cases} 0, & x \leq a^L; \\ \dfrac{x - a^L}{a^N - a^L}, & a^L \leq x \leq a^N; \\ \dfrac{a^R - x}{a^R - a^N}, & a^N \leq x \leq a^R; \\ 0, & x \geq a^R. \end{cases}$$

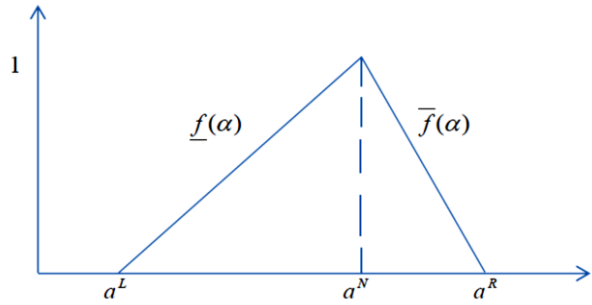

Fig. 1 Triangular Fuzzy Number (TFN)

The triangular fuzzy number $\tilde{A} = [a^L, a^N, a^R]$ may be transformed into interval form by using $\alpha$-cut in the following form

$$\tilde{A} = [a^L, a^N, a^R] = [a^L + (a^N - a^L)\alpha, a^R - (a^R - a^N)\alpha].$$

**Arithmetic**

Let us consider two fuzzy numbers $\tilde{x} = [\underline{x}(\alpha), \bar{x}(\alpha)]$ and $\tilde{y} = [\underline{y}(\alpha), \bar{y}(\alpha)]$ and a scalar $k$ then

a) $\tilde{x} = \tilde{y}$ if and only if $\underline{x}(\alpha) = \underline{y}(\alpha)$ and $\bar{x}(\alpha) = \bar{y}(\alpha)$.

b) $\tilde{x} + \tilde{y} = [\underline{x}(\alpha) + \underline{y}(\alpha), \bar{x}(\alpha) + \bar{y}(\alpha)]$.

c) $k\tilde{x} = \begin{cases} [k\underline{x}(\alpha), k\bar{x}(\alpha)], & k \geq 0, \\ [k\bar{x}(\alpha), k\underline{x}(\alpha)], & k < 0. \end{cases}$

If the fuzzy numbers are taken in order pair form as discussed earlier using limit method [10], the arithmetic rules are defined as

1. $[\underline{x}(\alpha), \bar{x}(\alpha)] + [\underline{y}(\alpha), \bar{y}(\alpha)]$
    $= [\min \{ \lim_{t \to \infty} m_1 + \lim_{t \to \infty} m_2, \lim_{t \to 1} m_1 + \lim_{t \to 1} m_2 \}, \max \{ \lim_{t \to \infty} m_1 + \lim_{t \to \infty} m_2, \lim_{t \to 1} m_1 + \lim_{t \to 1} m_2 \}]$

2. $[\underline{x}(\alpha), \bar{x}(\alpha)] - [\underline{y}(\alpha), \bar{y}(\alpha)]$
    $= [\min \{ \lim_{t \to \infty} m_1 - \lim_{t \to 1} m_2, \lim_{t \to 1} m_1 - \lim_{t \to \infty} m_2 \}, \max \{ \lim_{t \to \infty} m_1 - \lim_{t \to 1} m_2, \lim_{t \to 1} m_1 - \lim_{t \to \infty} m_2 \}]$

3. $[\underline{x}(\alpha), \bar{x}(\alpha)] \times [\underline{y}(\alpha), \bar{y}(\alpha)]$
    $= [\min \{ \lim_{t \to \infty} m_1 \times \lim_{t \to \infty} m_2, \lim_{t \to 1} m_1 \times \lim_{t \to 1} m_2 \}, \max \{ \lim_{t \to \infty} m_1 \times \lim_{t \to \infty} m_2, \lim_{t \to 1} m_1 \times \lim_{t \to 1} m_2 \}]$

4. $[\underline{x}(\alpha), \bar{x}(\alpha)] \div [\underline{y}(\alpha), \bar{y}(\alpha)]$
    $= [\min \{ \lim_{t \to \infty} m_1 \div \lim_{t \to 1} m_2, \lim_{t \to 1} m_1 \div \lim_{t \to \infty} m_2 \}, \max \{ \lim_{t \to \infty} m_1 \div \lim_{t \to 1} m_2, \lim_{t \to 1} m_1 \div \lim_{t \to \infty} m_2 \}]$

where for any arbitrary interval

$$[\underline{f}(\alpha), \overline{f}(\alpha)] = \left\{ \underline{f}(\alpha) + \frac{\overline{f}(\alpha) - \underline{f}(\alpha)}{t} = m \middle| \underline{f}(\alpha) \le m \le \overline{f}(\alpha), t \in [1, \infty) \right\}.$$

Let us consider $\tilde{A}(x)$ and $\tilde{B}(x)$ be two fuzzy numbers then it may be transformed into the following $\alpha$ cut form

$$\tilde{A}^{\alpha} = a^{\alpha} \text{ and } \tilde{B}^{\alpha} = b^{\alpha} \tag{10}$$

where

$$a^{\alpha} = [a_l^{\alpha}, a_r^{\alpha}] \text{ and } b^{\alpha} = [b_l^{\alpha}, b_r^{\alpha}]$$

$$a^{\alpha} = a_r^{\alpha} - \frac{a_r^{\alpha} - a_l^{\alpha}}{t} \text{ and } b^{\alpha} = b_r^{\alpha} - \frac{b_r^{\alpha} - b_l^{\alpha}}{t}, \ t \in [1, \infty).$$

Now consider an interval point $([a^L, a^R], [b^L, b^R])$ in two dimension plane and the pictorial representation of this number is shown in Fig. 2.

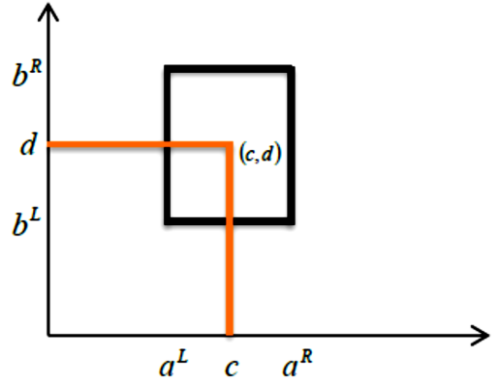

Fig. 2. Graphical representation of two-dimensional interval point

Here $[a^L, a^R]$ and $[b^L, b^R]$ may be represented in crisp form as defined in Eq. (1), $c$ and $d$ are the mid values of the intervals respectively. We get a set of crisp real points in the closed region i.e. $\{(a,b) \mid a = a^R - \frac{a^R - a^L}{t_1}, b = b^R - \frac{b^R - b^L}{t_2}; t_1, t_2 \in [1, \infty)\}$.

The above discussed fuzzy arithmetic is now used as a tool to solve FSDE. In the following sections the proposed technique is applied for both the exact and numerical methods to find the uncertain solutions of FSDE.

## 6. Solution of Fuzzy Stochastic Differential Equations (FSDE)

Let us consider a SDE with fuzzy parameters then Eq. (1) may be written as

$$d[\underline{X}(\alpha), \overline{X}(\alpha)] = [\underline{a}(\alpha), \overline{a}(\alpha)]dt + [\underline{b}(\alpha), \overline{b}(\alpha)]dW_t \tag{11}$$

Now Eq. (11) is solved by exact and numerical methods respectively.
Using limit method [2], the FSDE (11) in modified crisp form may be represented as follows

$$d\{\lim_{s \to \infty} X(\alpha), \lim_{s \to 1} X(\alpha)\} = \{\lim_{s \to \infty} a(\alpha), \lim_{s \to 1} a(\alpha)\}dt + \{\lim_{s \to \infty} b(\alpha), \lim_{s \to 1} b(\alpha)\}dW_t \tag{12}$$

where

$$X(\alpha) = \underline{X}(\alpha) + \frac{\overline{X}(\alpha) - \underline{X}(\alpha)}{s}, \ a(\alpha) = \underline{a}(\alpha) + \frac{\overline{a}(\alpha) - \underline{a}(\alpha)}{s} \text{ and } b(\alpha) = \underline{b}(\alpha) + \frac{\overline{b}(\alpha) - \underline{b}(\alpha)}{s}.$$

Initially for the exact case we take the crisp representation of $X(\alpha), a(\alpha), b(\alpha)$ and use Ito integral to solve the problem.

Now if we apply the above discussed fuzzy concept for Euler-Maruyam method, then Eq. (8) may be represented in the following way

$$X_0(\alpha) = w_0(\alpha)$$
$$w_{i+1}(\alpha) = w_i(\alpha) + a(t_i, w_i, \alpha)\Delta t_{i+1} + b(t_i, w_i, \alpha)\Delta W_i \quad (13)$$

where

$$X_0(\alpha) = \underline{X_0}(\alpha) + \frac{\overline{X_0}(\alpha) - \underline{X_0}(\alpha)}{s}, \quad w_0(\alpha) = \underline{w_0}(\alpha) + \frac{\overline{w_0}(\alpha) - \underline{w_0}(\alpha)}{s},$$

$$w_{i+1}(\alpha) = \underline{w_{i+1}}(\alpha) + \frac{\overline{w_{i+1}}(\alpha) - \underline{w_{i+1}}(\alpha)}{s}, \quad a(t_i, w_i, \alpha) = \underline{a}(t_i, w_i, \alpha) + \frac{\overline{a}(t_i, w_i, \alpha) - \underline{a}(t_i, w_i, \alpha)}{s}$$

and $b(t_i, w_i, \alpha) = \underline{b}(t_i, w_i, \alpha) + \frac{\overline{b}(t_i, w_i, \alpha) - \underline{b}(t_i, w_i, \alpha)}{s}$.

Applying $\lim_{s \to \infty}$ and $\lim_{s \to 1}$ on the solution we get the left and right bound. Whereas, we obtain various solution set by considering different values of membership function $\alpha \in [0,1]$. It is noticed that sometimes we get weak solutions i.e. the left and right bound solutions overlaps or intersect each other and this occurs due to the randomness of the system. This may easily be observed from the following example problems.

## 7. Example problems

In this section we have considered two example problems and taken parameters as fuzzy. Initially the problem is studied for crisp parameters for both the exact and numerical methods and then the problem is carrying out for fuzzy parameters.

Example 2
Consider Black Scholes stochastic differential equation [1].
The crisp Euler-Maruyama approximation for this SDE is as follows

$$w_0 = X_0$$
$$w_{i+1} = w_i + \mu w_i \Delta t_i + \sigma w_i \Delta W_i \quad (14)$$

where, the values of involved parameters are given in Table 1.

Table 1. Crisp and fuzzy values of the involved parameters.

| Parameters | crisp | TFN |
|---|---|---|
| $\mu$ | 0.75 | [0.65, 0.75, 0.85] |
| $\sigma$ | 0.30 | [0.25, 0.30, 0.35] |

Initially the Black Scholes SDE is solved for crisp parameter and then fuzzy parameters are considered for investigation. Here we compute a discretized Brownian path over [0, 1] with $\delta t = 2^{-8}$ and the obtained solution is plotted with a solid magenta line in Fig. 3. We then apply Euler Maruyama (EM) method using a step size $\Delta t = R\delta t$, with $R=4$ and obtain the solution which is presented in Fig 3 with blue line. Further it is seen that taking smaller value of $R$ for 4, 3 and 2 we get the endpoint errors 0.0442, 0.0216 and 0.0101 respectively.

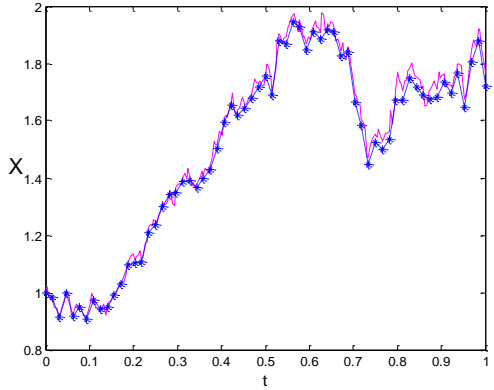
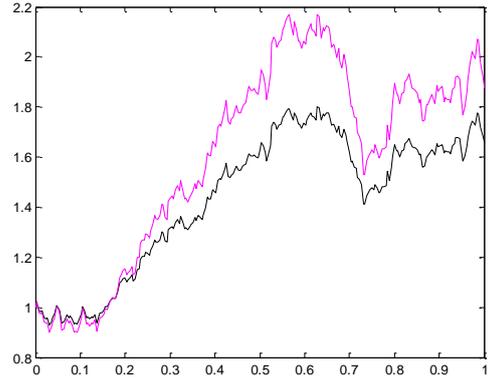

Fig 3. Solution of Black Scholes SDE when parameters are crisp

Fig 4. Exact solution of Black Scholes SDE when parameters are fuzzy.

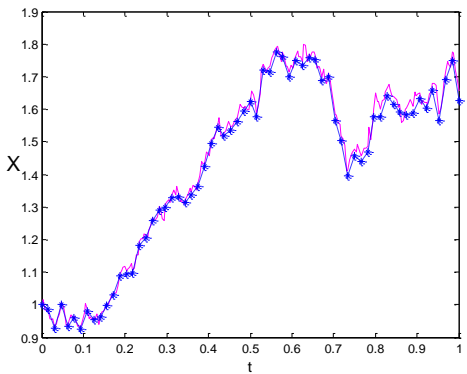
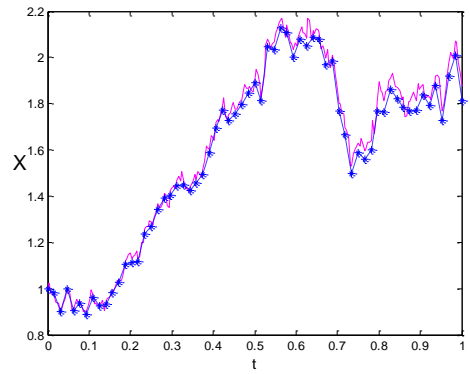

Fig 5. Crisp Euler Maruyama Solution of Black Scholes SDE with the left bound solutions.

Fig 6. Crisp Euler Maruyama Solution of Black Scholes SDE with the right bound solutions.

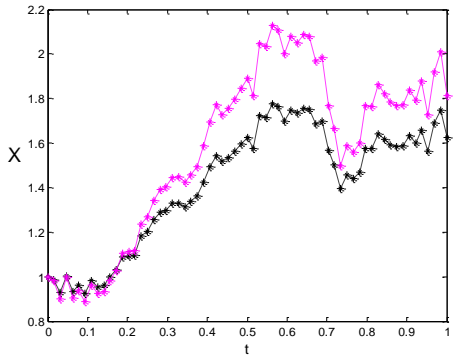
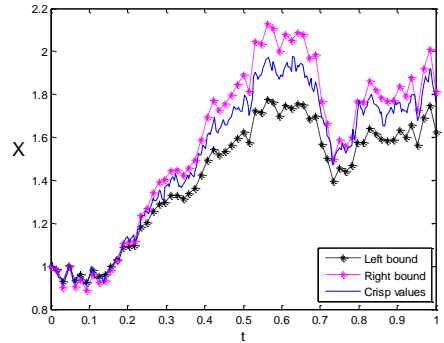

Fig 7. Euler Maruyama solution of Black Scholes SDE when parameters are fuzzy.

Fig 8. Euler Maruyama solution of Black Scholes SDE when parameters are fuzzy with the exact solution.

Now drift ($\mu$) and diffusion ($\sigma$) coefficients are taken as TFN which are given in Table 1. The exact method is used to obtain the solution which is depicted in Fig 4. Here the black and magenta solid line represents the left and right bound of the uncertainty. Next the left and right values of the uncertainty are plotted with the exact solution in Fig 5 and 6 respectively. Then EM method is used to solve the uncertain SDE and results are graphically depicted in Fig 7, where black and magenta line represents the left and right bound of the uncertain solutions. The region covered in between the left and right bound is the uncertain solution set of the Black Scholes SDE. In Fig 8, we have given the left and right bound is the uncertain solution of the Black Scholes SDE along with the crisp solution and we found that the exact solution lies within the region covered by the left and right solutions. Further it is found that there are some problems where it may be difficult to find the exact solution and in this case we take the help of Euler-Maruyama method which is discussed in the next example.

Example 3

The SDE of Langevin equation is
$$dX(t) = -\mu X(t)dt + \sigma dW_t \tag{15}$$
where $\mu$ and $\sigma$ are positive constants.

The Euler Maruyama approximation for Eq. (15) is as follows
$$\begin{aligned} w_0 &= X_0 \\ w_{i+1} &= w_i - \mu w_i \Delta t_i + \sigma \Delta W_i \end{aligned} \tag{16}$$

The values for used parameters for Eq. (12) are given in the following Table 2.

Table 2. Crisp and fuzzy values of the used parameters.

| Parameters | crisp | TFN |
|---|---|---|
| $\mu$ | 10 | [8, 10, 12] |
| $\sigma$ | 1 | [0.5, 1, 1.5] |

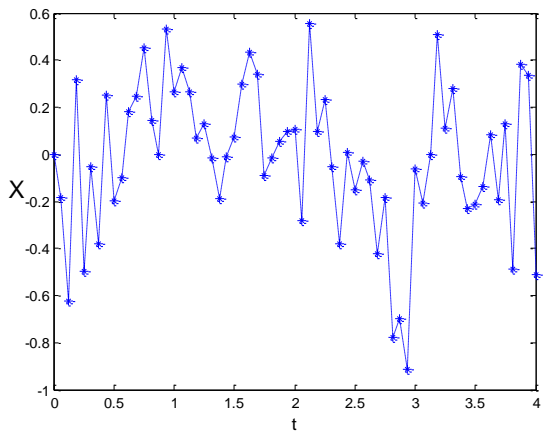
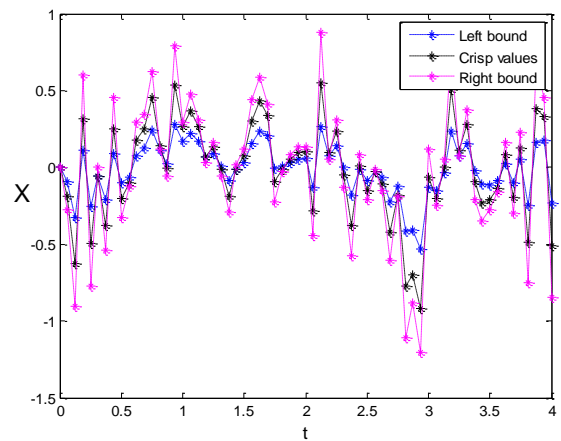

Fig 9. Euler Maruyama solution of Langevin SDE when parameters are fuzzy.

Fig 10. Euler Maruyama solution of Langevin SDE when parameters are fuzzy with the crisp solution.

In Fig 9 we have given a plot for the solution of Langevin SDE when parameters are crisp. Whereas, the solution for Langevin SDE is presented in Fig 10 where there parameters are taken as fuzzy. The left and right bound solutions are shown in blue and magenta colour respectively.

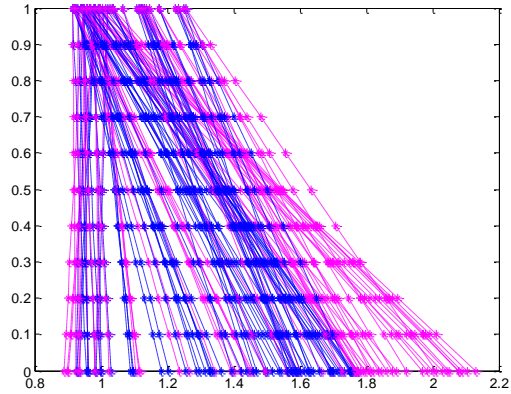 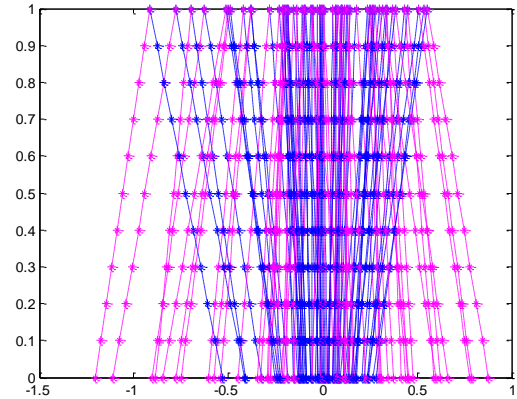

Fig 11. Fuzzy plot of Euler Maruyama solution of Black Scholes SDE when parameters are TFN (Example 2).

Fig 12. Fuzzy plot of Euler Maruyama solution of Langevin SDE when parameters are TFN (Example 3).

For better visualization of uncertain distribution of Euler Maruyama approximation results, fuzzy plots are represented in Fig 11 and 12 for example 2 and 3 respectively. Here the blue and magenta lines represent the left and right solutions and it is seen that uncertain widths are randomly varying.

It may be noted that if the uncertainty of the parameter changes, the uncertain width of the solution sets vary accordingly.

## 8. Conclusions

In this study we have combined the stochastic and fuzziness to model the problem and then proposed methods to handle them. Using $\alpha$-cut technique limit method is used to compute fuzzy arithmetic which is implemented to solve fuzzy stochastic differential equations both analytically and numerically. Then the hybrid method is demonstrated using example problems. Next the solution of uncertain width along with the left and right bound of the fuzzy stochastic differential equation have been investigated. Finally the TFN of the solution set is shown graphically for better visualization.

**Acknowledgements**

The authors would like to thank BRNS (*Board of Research in Nuclear Sciences*), Department of Atomic Energy, (DAE), Govt. of India for providing fund to do this work.